\documentclass[12pt]{article}
\usepackage{a4}
\usepackage{latexsym,amsmath,amssymb,amscd}
\usepackage{paralist}
\usepackage{amsthm}
\theoremstyle{plain}
    \newtheorem{theorem}                    {Theorem}       [section]
    \newtheorem{lemma}      [theorem]       {Lemma}
    \newtheorem{corollary}  [theorem]       {Corollary}
    \newtheorem{proposition}[theorem]       {Proposition}
    \newtheorem{definition} [theorem]       {Definition}
    \newtheorem{conjecture} [theorem]       {Conjecture}

    \newtheorem{question}   [theorem]       {Question}

\newtheorem{example}[theorem]{Example}

\newcommand{\alg}{{\text{alg}}}
\newcommand{\num}{{\text{num}}}
\newcommand{\homo}{{\text{hom}}}
\newcommand{\Hom}{\operatorname{Hom}}

\newcommand{\Gal}{\operatorname{Gal}}
\newcommand{\Pic}{\operatorname{Pic}}

\newcommand{\Tor}{{\operatorname{Tor}}}
\newcommand{\Cotor}{{\operatorname{Cotor}}}
\newcommand{\Spec}{\operatorname{Spec}}

\newcommand{\im}{\operatorname{im}}

\newcommand{\Br}{\operatorname{Br}}

\newcommand{\CH}{\operatorname{CH}}
\renewcommand{\lim}{\operatorname{lim}}
\newcommand{\colim}{\operatornamewithlimits{colim}}
\newcommand{\rk}{\operatorname{rank}}

\newcommand{\NS}{\operatorname{NS}}
\newcommand{\f}{{\mathcal F}}
\newcommand{\N}{{\mathbb N}}
\newcommand{\Z}{{{\mathbb Z}}}
\newcommand{\Q}{{{\mathbb Q}}}

\newcommand{\F}{{{\mathbb F}}}
\newcommand{\G}{{{\mathbb G}}}

\newcommand{\A}{{\mathcal A}}

\newcommand{\M}{{\mathcal M}}

\newcommand{\et}{{\text{\rm et}}}

\newcommand{\Div}{\operatorname{Div}}

\renewcommand{\div}{\operatorname{div}}
\newcommand{\corank}{\operatorname{corank}}
\renewcommand{\O}{{\cal O}}

\newcommand{\dash}{\text{-}}
\newcommand{\Sh}{\text{Sh}}

\makeatletter
\let\@fnsymbol\@arabic
\makeatother
\title{\bf On the structure of \'etale motivic cohomology.}
\author{Thomas H.~Geisser }

\begin{document}
\maketitle
\centerline{Dedicated to C.\ Weibel on his 65th birthday}

\begin{abstract}
We discuss the structure of integral \'etale motivic cohomology groups of
smooth and projective schemes over algebraically closed fields, 
finite fields, local fields, and arithmetic schemes. 
\end{abstract}

\section{Introduction}
Let $B$ be the spectrum of a Dedekind ring or field, $X$ a smooth and
projective variety over $B$, and $\Z(n)$ Bloch's cycle complex. The goal of
this paper is to find structure results for
the integral \'etale motivic cohomology groups $H^i_\et(X,\Z(n))$ of $X$. 
Our first result concerns algebraically closed base fields: 

\begin{theorem}\label{maintt}
Let $B=\Spec k$ be the spectrum of an algebraically closed field of 
characteristic $p\geq 0$, $i\not=2n$, and $d=\dim X$. Then we have an isomorphism 
$$H^i_\et(X,\Z(n))\cong D^i(n)\oplus C^i(n)\oplus P^i(n),$$
where
\begin{enumerate}
\item $D^i(n)$ is uniquely divisible, and vanishes for $i>\min\{2n,n+d\}$,
\item the prime-to-$p$ torsion part $C^i(n)$ is isomorphic to
$\Q/\Z[\frac{1}{p}]^r\oplus F^i(n)$, 
invariant under extensions of algebraically closed fields, 
$F^i(n)=\prod_{l\not=p}\Tor H^i_\et(X,\Z_l)$ is a finite group independent of $n$, 
$r\in \N_0$ is independent of $n$ if $n\not=\frac{i-1}{2}$, and
$C^i(n)$ vanishes unless $1\leq i\leq 2d+1$,
\item the $p$-primary torsion part 
$P^i(n)$ is isomorphic to $(\Q_p/\Z_p)^s\oplus P^i_c(n)$, with
$P^i_c(n)$ an extension of a finite group $F^i_p(n)$ by a 
finitely generated torsion $W(k)$-module $U^i(n)$, and 
$P^i(n)=0$ unless $n+1\leq i\leq n+d+1$ and $ 0\leq n\leq d$.
\end{enumerate}
\end{theorem}

We have the following duality statements:

\begin{proposition} \label{dualprop}
Under the hypothesis of the theorem, 
\begin{enumerate}
\item $F^i(n)$ and $F_p^i(n)$ are Pontrjagin dual to $F^{2d+1-i}(d-n)$
and $F_p^{2d+1-i}(d-n)$, respectively.
\item If $n,u\geq 0, i\not=2u+1$, and $ 2d+2-i\not=2n+1$, then the Tate module of 
$C^{2d+2-i}(n)$ is Pontrjagin dual to $C^i(u)/F^i(u)$.
\item $U^i(n)\cong \Hom_{W(k)}(U^{2d+2-i}(d-n),CW(k))$, for 
$CW(k)$ the co-Witt vectors.
 \end{enumerate}
\end{proposition}

From the Rost-Voevodsky theorem it follows that 
the motivic cohomology groups $H^i_\M(X,\Z(n))$ satisfy
the conclusion of the previous theorem
for $i\leq n+1$ and for $n\geq \dim X$ as well.
We mention some results and examples on $H^{2n}_\et(X,\Z(n))$ as well.

If the base field is finite, then a conjecture of Lichtenbaum gives a
precise conjectural description of 
the structure of \'etale motivic cohomology \cite{lichtenbaumMC}, and 
the conjecture is equivalent to other deep and well-known conjectures. 
Over local fields, our main result is:

\begin{theorem}\label{localstructure}
Let $k$ be a $p$-adic field of residue characteristic $p$. 
Then $H^i_\et(X,\Z(n))$ is
the direct sum of a finite group and a group which is uniquely $l$-divisible 
for all $l\not=p$, if either
$X$ has good reduction and $i\not\in \{2n-1, 2n,2n+1, 2n+2\}$, or if 
$i\not \in \{n,\ldots , n+d+2\}$.
\end{theorem}

We also give a list of examples showing that the bounds are sharp,
and ask some more precise questions on the structure of the groups. 

Finally, let $C$ be spectrum of the ring of integers of a number field or a smooth
and proper curve over a finite field. In analogy to the situation over
finite fields, Lichtenbaum conjectures that if $X$ is regular, and proper
over $C$, then the groups  
$H^{i}_\et(X,\Z(n))$ are finitely generated for $i\leq 2n$,
finite for $i=2n+1$ and of cofinite type for $i\geq 2n+2$. 
If $B=C-S$ is the complement of a finite, non-empty set $S$ of places, and $X$
is smooth and proper over $B$, then we expect $H^{i}_\et(X,\Z(n))$
to be the direct sum of a finitely generated group and a group of
cofinite type, and we raise the question if the group
$$\Sh^{i,n}(X)= \ker  H^i_\et(X,\Z(n)) \to 
\bigoplus_{v\in S} H^i(X\times_BK_v,\Z(n)),$$
where $K_v$ is the completion of the function field of $B$
at the place $v$, 
is finite if $S$ contains at least one finite place. 
We show that the answer is (trivially) affirmative
for some small values of $i$ and $n$, and prove

\begin{theorem}
If $f:X\to B$ is smooth and proper,
then $\Sh^{3,1}(X)$ is finite if and only if 
the Tate-Shafarevich group $\Sh(\Pic^{0,red}_{X_\eta/\eta})$ of
the Picard variety of the generic fiber is finite.
\end{theorem}

\medskip

{\bf Notation:} Throughout the paper, $n$ will be a non-negative integer.
We denote $\Z(n)$ Bloch's motivic complex of cycles of codimension $n$, 
a complex of \'etale sheaves on the 
category $Sm/B$ of smooth schemes over $B$, \cite[Thm.1.17]{kahnglr},
by $\Z(n)$. 

For an abelian group $A$ we denote by ${}_m A$ its $m$-torsion,
by $A\{l\}=\colim_{l^r} {}_{l^r}A$ its subgroup of
$l$-power torsion elements, by $A^*$ its Pontrjagin dual
$\Hom(A,\Q/\Z)$, by $A^\wedge=\lim A/m$ its completion,
by $A^{\wedge l}=\lim A/l^r$ its $l$-adic completion, by 
$TA=\lim {}_{m} A$ its Tate module, and by 
$T_lA=\lim_r {}_{l^r} A$ its $l$-adic Tate module.
The subgroup of $l$-divisible elements is denoted by $l\dash \div A$
and the largest $l$-divisible subgroup by $l\dash \Div A$.

We will call an $l$-power torsion group
of cofinite type, if it is of the form $(\Q_l/\Z_l)^r\oplus F$
for a finite group $F$, and we call a torsion group of cofinite type 
if it is of the form $\Q/\Z[\frac{1}{p}]^r\oplus \Q_p/\Z_p^s$
(where $p$ is the characteristic of the base field in case of algebraically
closed and finite fields, and the residue characteristic on case of local fields).

\section{Algebraically closed base fields}
Assume that $k$ is an algebraically
closed field of characteristic $p\geq 0$. 

\begin{proposition}\label{cdrk}
If $l\not=p$ is a prime number, then 
the $l$-adic cohomology 
groups $H^i_\et(X,\Z_l(n))$ are finitely generated $\Z_l$-modules, 
of rank independent of $l\not=p$, and torsion free for
almost all $l$. The groups $H^i_\et(X,\Q_l/\Z_l(n))$
and $H^i_\et(X,\Z(n))\{l\}$ are of cofinite type, of corank independent
of $l$, and cofree for almost all $l$.
\end{proposition}

\proof  
By SGA 4.5, the groups $H^i_\et(X,\Z/l^r(n))$ are finite, hence
$H^i_\et(X,\Z_l(n))= \lim H^i_\et(X,\Z/l^r(n))$ is a compact 
$\Z_l$-module, and
hence finitely generated. By Gabber \cite{gabber}, 
its $l$-torsion vanishes for almost all $l$.
The rank does not depend on $l$ by comparing with the Betti-number
in characteristic $0$, and by the Weil-conjectures in characteristic $p$.
The statements for $\Q_l/\Z_l$-coefficients follows by taking the colimit,
and the statement about torsion 
follows from the surjection 
$ H^i_\et(X,\Q_l/\Z_l(n))\twoheadrightarrow H^{i+1}_\et(X,\Z(n))\{l\}$.
\endproof

For the $p$-part, we have the following:

\begin{proposition}\label{ccrk}
The groups $H^i_\et(X,\Z_p(n))$ are the direct sum of a finitely
generated free $\Z_p$-module and an extensions of a finite group 
by a finitely generated torsion $W(k)$-module. 
The groups $H^i_\et(X,\Q_p/\Z_p(n))$ and 
$H^i_\et(X,\Z(n))\{p\}$ 
are the direct sum of a group of the form $(\Q_p/\Z_p)^r$, and
an extension of a finite group by a finitely generated torsion
$W(k)$-module.
\end{proposition}

\proof
We have $H^i_\et(X,\Z/p^r(n))\cong H^{i-n}_\et(X,\nu_r^n)$
by \cite{ichmarc}. By Milne \cite[Lemma 1.8]{milneppart}, 
the sheaf associated to the
presheaf $H^j(X\times -,\nu_r^n)$ on the category $Pf/k$ of perfect schemes
over $k$ with the \'etale topology is represented by a commutative perfect
group scheme over $k$, corresponding to an extension of an \'etale
group scheme by a unipotent commutative quasi-algebraic group,
whose limit is finite dimensional by loc.cit. Prop. 3.1.
Taking global sections over $k$, and then the
limit, the results follows. The statement
about torsion coefficients follows by taking the colimit, and 
the statement on torsion follows from the short exact sequence
$$0\to H^i_\et(X,\Z(n))\otimes \Q_p/\Z_p\to 
 H^i_\et(X,\Q_p/\Z_p(n))\twoheadrightarrow H^{i+1}_\et(X,\Z(n))\{p\}\to 0.$$
Indeed, since the left hand group is divisible, the cotorsion of
the two groups on the right is isomorphic.
\endproof

\begin{corollary}
The torsion subgroup $\Tor H^{i}_\et(X,\Z(n))$ is a direct summand
of $H^{i}_\et(X,\Z(n))$.
\end{corollary}

\proof
The proposition shows that $ \Tor H^{i}_\et(X,\Z(n))$ is the direct
sum of a divisible group and a group of finite exponent. The result
now follows from \cite[Thm. 16]{kap}.
\endproof

\begin{proposition}\label{suslin}
Let $F/k$ be an extension of algebraically closed fields. 
Then for $l\not=p$, the base change maps
\begin{align*}
H^i_\et(X,\Z(n))\{l\}&\stackrel{\sim}{\to} H^i_\et(X_F,\Z(n))\{l\}\\
H^i_\et(X,\Z(n))^{\wedge l}&\stackrel{\sim}{\to} 
H^i_\et(X_F,\Z(n))^{\wedge l}
\end{align*}
are isomorphisms. For $l=p$, then are injective with cokernel of bounded exponent.
Furthermore we have
$$H^i_\et(X,\Z(n))\otimes\Q/\Z\cong H^i_\et(X_F,\Z(n)\otimes\Q/\Z.$$
\end{proposition}

\proof
Consider the map of short exact coefficient sequences:
$$\begin{CD}
H^i_\et(X,\Z(n))/l^r@>>> H^i_\et(X,\Z/l^r (n))@>>> {}_{l^r}H^{i+1}_\et(X,\Z(n))\\
@VVV @VVV @VVV \\
H^i_\et(X_F,\Z(n))/l^r@>>> H^i_\et(X_F,\Z/l^r (n))@>>> 
{}_{l^r}H^{i+1}_\et(X_F,\Z(n))
\end{CD}$$
Since the middle map is an isomorphism by the smooth and proper base
change theorem for $l\not=p$, and is injective with a cokernel bounded
independently of $r$ for $l=p$ because
$H^i_\et(X,\Z/p^r(n))\cong H^{i-n}_\et(X,\nu_r^n)$ and
\cite[Lemma 1.8]{milneppart},
it suffices to show that the outer maps are injective. 
We write $F$ as a colimit of finitely 
generated $k$-algebras $A$, and note that $k\subseteq A$ is
split because $k$ is algebraically closed. Since 
\'etale cohomology commutes with limits with affine transition maps, 
the result follows. 
\endproof

Note that no information can be deduced about a group 
from knowing its completion and torsion subgroup, even in the absence
of uniquely divisible groups. 
For example, $G=\Z$ and $G=\oplus_l \Z_{(l)}$ both
satisfy $G^{\wedge l}\cong \Z_l$
and $G\otimes \Q_l/\Z_l\cong \Q_l/\Z_l$ for all $l$.

\medskip

\proof (Theorem \ref{maintt}).
It suffices to show that $H^i_\et(X,\Z(n))\otimes \Q/\Z=0$. Indeed, 
then the short exact sequence 
$$ 0\to \Tor H^i_\et(X,\Z(n))\to H^i_\et(X,\Z(n))\to H^i_\et(X,\Z(n))\otimes \Q
\to 0$$
shows that $H^i_\et(X,\Z(n))$ modulo its torsion subgroup is uniquely 
divisible, and the structure of the torsion subgroup is given 
in Propositions \ref{cdrk} and \ref{ccrk}.

To show $H^i_\et(X,\Z(n))\otimes \Q/\Z=0$, we take a (normal, reduced) scheme $S$
of finite type over $\Z$ and a scheme $\mathcal X$ over $S$ 
such that $X$ is the base change to $k$ 
of the generic fiber $\mathcal X_\eta$.
From Proposition \ref{suslin} we know that 
$H^i_\et(\mathcal X_{\bar\eta},\Z(n))\otimes\Q/\Z\cong H^i_\et(X,\Z(n))\otimes
\Q/\Z$, hence it suffices to consider the case $k=\bar\eta$.
Let $N$ be $H^i_\et(X,\Z(n))$ modulo its torsion subgroup, it suffices to show
that  $N$ is divisible. From the
exact sequence $0\to l\dash \div N\to N\to N^{\wedge l}$ we see that it is enough 
to show that for all $l$ the map 
$$N\to N^{\wedge l} \subseteq H^i_\et(X,\Z_l(n))/\Tor 
\subseteq H^i_\et(X,\Q_l(n))$$
is the zero map. Since the smooth locus of $\mathcal X$ is
open (and non-empty because $X$ is smooth), we can find a point 
$\Spec \F$ of $S$ with finite residue field such that the fiber
$\mathcal X_\F$ is smooth. If $l\not=p$, we have an isomorphism  
$H^i_\et(X_{\bar\F},\Q_l(n)) \cong H^i_\et(X_{\bar \eta},\Q_l(n))$, 
compatible with the action of $\Gal(\eta)\to \Gal(\F)$,
by the smooth and proper base change theorem. For $l=p$, 
we can assume that $S$ is of finite type over a finite
field contained in $\F$, and the same statement holds by Gros-Suwa 
\cite[Thm.\ II 2.1]{grossuwa}. 
Now since
$$N\subseteq H^i_\et(X,\Q(n))\cong H^i_\M(X,\Q(n))= 
\colim_{L\supseteq \eta\; \text{finite}}  H^i_\M(X_L,\Q(n)),$$
every element of $N$ has finite orbit under $\Gal(\eta)$.
In particular, the arithmetic Frobenius $\varphi$ has finite order
acting on any element of $N$.
On the other hand, by the Weil conjectures and its $p$-adic version, 
the eigenvalues of $\varphi$ acting
on $H^i_\et(X_{\bar\F},\Q_l(n))$ have absolute value $p^{n-\frac{i}{2}}$,
hence there is no element of finite order. 
\endproof

\proof (Proposition \ref{dualprop}). 
Since $T_l H^{i+1}_\et(X,\Z(n))$ is torsion free, 
we get from the coefficient sequence that 
$$F^i(n)\oplus P^i_c(n) \cong \Tor H^i_\et(X,\Z(n))^\wedge \cong \Tor  \prod_{l}
H^i_\et(X,\Z_l(n)).$$
If $l\not=p$, then $F^i(n)$ does not depend on $n$ because 
$H^i_\et(X,\Z_l) \cong H^i_\et(X,\Z_l(n))$ for every $n$.
Taking the inverse limit of Poincar\'e-duality we obtain
$$H^i_\et(X,\Z_l(n))\cong H^{2d-i}_\et(X,\Q_l/\Z_l(d-n))^*.$$
Since $\Hom(D,\Q/\Z)$ is uniquely divisible for a divisible group $D$, the map 
$$ H^{2d-i+1}_\et(X,\Z_l(d-n))^*\to 
H^{2d-i}_\et(X,\Q_l/\Z_l(d-n))^*$$
induces an isomorphism 
$$\big( \Tor H^{2d-i+1}_\et(X,\Z_l(d-n))\big) ^*\cong  
\big( \Cotor\; H^{2d-i}_\et(X,\Q_l/\Z_l(d-n))\big) ^*.$$

For the $p$-part, we use Milne's duality \cite[Thm.1.11]{milneppart}, see also 
\cite[Cor. 3.26]{grossuwa}. 
Let $\underline U^i(X,\nu_r(n))$ be the unipotent part  and 
$\underline D^i(X,\nu_r(n))$ be the \'etale quotient 
of the finite dimensional
pro-group scheme $\underline H^{i-n}(X\times -,\nu_r(n))$ on $Pf/k$,
see also \cite[Cor. 3.25]{grossuwa}.
Taking $k$-rational points we obtain $H^{i}_\et(X,\Z/p^r(n))$, and 
since $\underline U^i(X,\nu_\cdot(n))$ is finite dimensional by 
loc.cit. Prop. 3.1, we get
$U^i(n)\cong \colim_r \underline U^{i-n-1}(X,\nu_r(n))(k)$,
and $F^i_p(n)$ is isomorphic both to the cotorsion of 
$\colim_r \underline D^{i-n-1}(X,\nu_r(n))(k)$ and to the torsion 
in $\lim_r \underline D^{i-n}(X,\nu_r(n))(k)$.
Hence $D^i(n)\cong \Hom(D^{2d+1-i}(d-n),\Q_p/\Z_p)$, which gives
duality for $F^i_p(n)$ exactly as above, and 
$U^i(n)\cong \Hom_{W(k)}(U^{2d+2-i}(d-n),CW(k))$.
\endproof
 
\subsubsection*{Remarks}
\begin{enumerate}
\item 
The method of Colliot-Th\'el\`ene and Raskind \cite[Th.\ 1.8, 2.2]{ctras}
yields similar results away from the characteristic $p$.
\item It would be interesting to write down the duality pairing
between $U^i(n)$ and $U^{2d+2-i}(d-n)$ in terms of the motivic
cohomology groups directly.
\item
The Beilinson-Soul\'e vanishing conjecture is equivalent to the
vanishing of $ D^i(n)=0$ for $i<0$.
\item
If the base field has characteristic $p>0$ 
and if we assume Parshin's conjecture, then $D^i(n)=0$ for $i<n$ \cite{ichtate}. 
A more careful analysis shows that 
if the base-field has transcendence degree $r$ over the finite prime field, 
then $H^i_\et(X,\Q(n))=0$ for $i<\max\{n,2n-r\}$ under Parshin's conjecture.
\item There is no non-degenerate pairing between $D^i(n)$ and some $
D^j(d-n)$, because for the algebraic closure of the rationals, 
$D^1(n)$ is infinite dimensional for every $n>0$.
\end{enumerate}

\subsubsection*{Examples}

\begin{enumerate}
\item \label{fgleichnull}
We have $F^0(n)=F^1(n)=F^{2d}(n)=0$, $F^2(n)$ is the 
prime to $p$-torsion of the Neron-Severi group, and $F^3(n)$ is
the prime to $p$-cotorsion of the Brauer group. 
\item If $X$ is a supersingular abelian surface or K3 surface, then
$H^3_\et(X,\Z(1))$ has a direct summand isomorphic  to $k$, dual to itself.
\item We have 
$$H^1_\et(\Q,\Z(1))\cong \Q^\times \cong \Z/2 \oplus \bigoplus_p \Z$$
is the Galois invariants of 
$$H^1_\et(\bar \Q,\Z(1))\cong 
 \bar\Q^\times\cong  \Q/\Z\oplus\bigoplus_{\mathfrak p}\Q ,$$
where $\mathfrak p$ runs over all finite places of $\bar\Q$. 
This can be explained by the Galois cohomology sequence associated
to the short exact sequence of Galois modules
$$ 0\to \Tor H^i_\et(\bar X,\Z(n)) \to H^i_\et(\bar X,\Z(n)) \to D^i(n) \to 0.$$
In particular, the
action of the Galois group is not compatible with the decomposition in 
Theorem \ref{maintt}. 
\end{enumerate}

\section{The \'etale Chow-group}
This section gives an overview over some known results in 
degree $2n$.

\subsection*{Equivalence relations}
Let $X$ be again a smooth and projective scheme over an algebraically
closed field.
The intersection and cup product gives us a diagram of pairings
\begin{equation}\label{ddd}
\begin{CD}
\CH^n(X) \times \CH^{d-n}(X) @>\cap >> \CH^d(X)@>>> \Z \\
@VVV @VVV @| \\
H^{2n}_\et(X,\Z(n))\times H^{2d-2n}_\et(X,\Z(d-n))@>\cup>> 
H^{2d}_\et(X,\Z(d))
@>>> \Z\\
@VVV @VVV @VVV \\
H^{2n}_\et(X,\Z_l(n))\times H^{2d-2n}_\et(X,\Z_l(d-n))@>\cup >> 
H^{2d}_\et(X,\Z_l(d))
@>>> \Z_l,
\end{CD}\end{equation}
and the kernels of these pairings are the subgroups of elements 
numerically equivalent to zero,
$\CH^n_\num(X)\subseteq \CH^n(X)$ and 
$H^{2n}_\num(X,\Z(n))\subseteq H^{2n}_\et(X,\Z(n))$. 
The lower pairing is non-degenerate modulo torsion. 

\begin{proposition}\label{modnum}
The map $$\CH^n(X)/\CH^n_\num(X)\to  
H^{2n}_\et(X,\Z(n))/H^{2n}_\num(X,\Z(n))$$
is an isomorphism of finitely generated free abelian groups. 
\end{proposition}

\proof
We write  $A^D=\Hom(A,\Z)$ for an abelian group $A$.
We have non-degenerate pairings on 
$N= \CH^n(X)/\CH^n_\num(X)$ 
to $\Z$, which shows that these groups are torsion free. 
On the other hand, since the pairings above are compatible, and 
the lower pairing is perfect after tensoring with $\Q_l$, 
$N\otimes\Q_l$ injects into the finite
dimensional $\Q_l$-vector space $H^{2n}_\et(X,\Q_l(n))$.
This implies that 
$N\otimes\Q$ is finite dimensional. Let $\{n_1, \ldots , n_r\}$ be
elements of $N$ which form a basis of $N\otimes \Q$, and $N_1\subseteq N$
be the submodule generated by $\{n_1, \ldots , n_r\}$. 
Since $N/N_1$ is torsion, we have $(N/N_1)^D=0$, hence the map 
$$M= \CH^{d-n}(X)/\CH^{d-n}_\num(X)\subseteq  N^D \to N_1^D$$ 
is injective. Thus $M$ injects into a finitely generated free abelian group,
hence is finitely generated free.
The same argument applies to 
$H^{2n}_\et(X,\Z(n))/ H^{2n}_\num(X,\Z(n))$.
To get the isomorphism, consider the commutative diagram
$$\begin{CD}
\CH^n(X)/\CH^n_\num(X) @>>>  H^{2n}_\et(X,\Z(n))/ H^{2n}_\num(X,\Z(n))\\
@| @| \\
\big(\CH^{d-n}(X)/ \CH^{d-n}_\num(X)\big)^D@<<< 
\big( H^{2d-2n}_\et(X,\Z(d-n))/ H^{2d-2n}_\num(X,\Z(d-n))\big)^D.
\end{CD}$$
Since $\CH^n(X)_\Q \cong H^{2n}_\et(X,\Q(n))$, the upper map is
rationally surjective, hence all groups are finitely generated
free abelian groups of the same rank, and it follows that all maps 
are isomorphisms.
\endproof

We define the groups $\CH^n_\homo(X)$ and $H^{2n}_\homo(X,\Z(n))$ 
of cycles homologically equivalent to zero to be the kernel of 
the composition, and of the second map in
$$ \CH^n(X)\to H^{2n}_\et(X,\Z(n)) \to \prod_l H^{2n}_\et(X,\Z_l(n)),$$
respectively.
Since the second map factors through the subgroup $\lim_m H^{2n}_\et(X,\Z(n))/m$,
$H^{2n}_\homo(X,\Z(n))$ is the group of 
divisible elements.

\begin{lemma}
The group $H^{2n}_\homo(X,\Z(n))$ is the maximal divisible subgroup of 
$H^{2n}_\et(X,\Z(n))$, and the map 
$$\CH^n(X)/ \CH^n_\homo(X)\to 
H^{2n}_\et(X,\Z(n))/H^{2n}_\homo(X,\Z(n))$$
is injective with torsion cokernel.
\end{lemma}

\proof
If an element in $\CH^n(X)$ maps to $H^{2n}_\homo(X,\Z(n))$,
then it vanishes in $l$-adic cohomology, hence is contained
in $\CH^n_\homo(X)$ which implies injectivity. The cokernel
is torsion because the Chow group and \'etale Chow group agree
rationally.
\endproof

It is conjectured that $\CH^n_\homo(X)\subseteq \CH^n_\num(X)$ agree
up to a torsion group.
The same argument as in the Lemma gives:

\begin{lemma}
The map 
$$\CH^n_\num(X)/ \CH^n_\homo(X)\to 
H^{2n}_\num(X,\Z(n))/H^{2n}_\homo(X,\Z(n))$$
is injective with torsion cokernel. 
\end{lemma}

A cycle is algebraically equivalent to zero if it lies in the 
image of some map
$$\begin{CD}
\CH^n(T\times X) @>t_1^*-t_0^*>> \CH^n(X)\\
@VVV @VVV \\
H^{2n}_\et(T\times X,\Z(n)) @>t_1^*-t_0^*>> H^{2n}_\et(X,\Z(n)),
\end{CD}$$
where $T$ is a smooth connected scheme $T$ 
(which we can assume to be a smooth
curve) with closed points $t_0,t_1\in T$.

We obtain compatible filtrations on Chow groups
and \'etale Chow groups 
$$\begin{CD} 0@.\subseteq \CH^n_\alg (X)@.\subseteq \CH^n_\homo(X) 
@.\subseteq \CH^n_\num(X) @.\subseteq \CH^n(X)\\
@. @VVV @VVV @VVV @VVV \\
0@.\subseteq H^{2n}_\alg(X,\Z(n)) @.\subseteq H^{2n}_\homo(X,\Z(n)) 
@.\subseteq H^{2n}_\num(X,\Z(n)) @.\subseteq H^{2n}(X,\Z(n)).
\end{CD}$$
The upper and lower groups agree rationally.

\subsection*{Examples}
\begin{enumerate}
\item
(Bloch and Esnault \cite{blochesnault}) There are 
$3$-dimensional complete
intersections $X$ over the algebraic closure of a number field
such that $\CH^2_\homo(X)$ is not $l$-divisible,
and $CH^2(X)\{l\}=0$ for some primes $l$.
\item
(Schoen 
\cite{schoenjag})
The triple self-product of an elliptic 
curve over an algebraically closed field of characteristic $0$
can have infinite $\CH^n(X)/l$ for $2\leq n\leq d-1$, in particular 
$\CH^n_\homo(X)/l$ is infinite. 
From the injectivity of 
$\CH^2_\homo(X)\subseteq \CH^2(X)\to H^{4}_\et(X,\Z(2))$ 
and finiteness of $H^{4}_\et(X,\Z/l(2))$
it follows that there is a subgroup of infinite rank
in $\CH^2(X)$ which becomes divisible in $H^{4}_\et(X,\Z(2))$.
\item
(Schoen 
\cite{schoenjams}) Over the algebraic 
closure of a finite field of characteristic larger than $2$, 
and assuming Tate's conjecture and semi-simplicity of the
Frobenius action, the map
$\CH^n_\homo(X) \to H^{2n-1}_\et(X,\Q_l/\Z_l(n))$ 
is surjective for almost all $l\not=p$ and $n$ if the dimension of $X$ is at
most $4$, and in general under some additional hypothesis.
Note that over the algebraic closure of a finite field,
$\CH^n_\homo(X)$ is conjectured to be a torsion group.
\item
(Schoen 
\cite{schoenma}) For any $n\geq 2$ there are varieties of
dimension $2n-1$ such that the Griffiths group
$Gr^n(X)=\CH^n_\homo(X)/\CH^n_\alg(X)$ has torsion elements.
\item
(Schoen 
\cite{schoenmrl}) 
There are smooth, projective varieties $V$, 
such that the corank of $\CH^r(X)\{l\}$ is infinite for all $r$ 
in the range $2 < r < \dim(X)$.
\end{enumerate}

\subsection*{Representability}
Recall that a group homomorphism from $\CH^n_\alg(X)$ to the 
$k$-rational points of an abelian
variety $A$ is regular, if for every pointed smooth connected
variety $t_0\in T$ and 
correspondence  $\Gamma\in \CH^n(T\times X)$, the composition
with the map 
$$T(k)\to \CH^n_\alg(X),
\qquad t\mapsto t^*\Gamma-t_0^*\Gamma$$
is the map induced on closed points by a morphism of schemes $T\to A$. 
We use the same concept for the \'etale Chow group:

\begin{definition}
A homomorphism from $H^{2n}_\alg(X,\Z(n))$ to the $k$-rational points of an abelian
variety $A$ is regular, if for every pointed smooth connected 
variety $t_0\in T$ and element $\Gamma\in H^{2n}_\alg(T\times X,\Z(n))$, 
the composition with
$$T(k)\to H^{2n}_\alg(X,\Z(n)), \qquad t\mapsto t^*\Gamma-t_0^*\Gamma$$
is the map induced on closed points by a morphism of varieties $T\to A$.
\end{definition}

\begin{theorem}\label{repr}
There is a universal object $\rho_n:H^{2n}_\alg(X,\Z(n))\to A_n$ 
for regular homomorphisms from $H^{2n}_\alg(X,\Z(n))$ to abelian varieties.
\end{theorem}

\proof
This follows by the argument of Serre \cite{serre} and H.\ Saito \cite{saito} 
because the dimension of surjective maps to abelian varieties is 
bounded.
\endproof

In \cite{murrebar}, \cite{murretor}, Murre studied the situation
for Chow groups, and he proved that a universal homomorphism to abelian
varieties exists for dimension $0$, and codimensions $1$ and $2$,
so also \cite{korita}.

\section{Other bases}
\subsection*{Finite fields}
We mention a conjecture on the structure of 
motivic cohomology groups over finite fields:

\begin{conjecture}\label{lconj} (Lichtenbaum \cite{lichtenbaumMC}) 
$$ H^i_\et(X,\Z(n))=
\begin{cases}
\text{finite} &i\not=2n, 2n+2, \\
\text{finitely generated} & i=2n, \\
\text{cofinite type} &i=2n+2.
\end{cases}$$
\end{conjecture}

We note that a small modification of \'etale motivic cohomology,
called Weil-\'etale cohomology, 
yields groups which are conjecturally finitely generated for all $i$,
\cite{lichtenbaumweil, ichweil}.

\begin{example}\label{ff} For $X$ a point we have
$$ H^i_\et(\F_q,\Z(n))= 
\begin{cases}
\Z & (i,n)=(0,0) \\
\Q/\Z  &  (i,n)=(2,0)\\
\Z/{q^n-1} & i=1, n>0\\
0 & \text{otherwise}
\end{cases}$$
For $n=0$ this is the calculation of Galois cohomology,
and for $n\geq 1$, the groups $H^i_\et(\F_q,\Z(n))$
are annihilated by $q^n-1$ by \cite[Thm. 4.6]{ichparshin}, hence they are isomorphic to the well-known groups $H^{i-1}_\et(\F_q,\Q/\Z(n))$.
\end{example}

\begin{proposition}\label{ccc}
Lichtenbaum's conjecture is true for curves. 
\end{proposition}

\proof
By Soul\'e \cite{soulecurves}, $H^i_\et(X,\Z(n))$ is torsion
unless $(i,n)=(0,0)$ or $(2,1)$, in which case the groups are $\Z$
and $\Z\oplus (finite)$, respectively. On the other hand, 
$H^i_\et(X,\Q/\Z(n))$ is of cofinite type, and finite unless 
$(i,n)= (0,0), (1,0), (2,1)$, or $ (3,1)$. It follows that
$\Tor H^i_\et(X,\Z(n))$ is finite unless $(i,n)= (1,0), (2,0), (3,1)$, or
$(4,1)$.
Finally, $H^1_\et(X,\Z(0))=0$ and finiteness of $H^3_\et(X,\Z(1))=\Br(X)$, follows
from Tate's theorem.
\endproof

\begin{proposition}
For $i<2n$, Conjecture \ref{lconj} is true up to a uniquely divisible group,
and true if and only Parshin's conjecture holds. For $i\geq 2n+2$, it
holds unconditionally, and for 
$i=2n+1$, it is equivalent to  Tate's conjecture
on the surjectivity of the cycle map.

The conjecture holds for all $X$ and $n$ in all degrees if and only if 
Tate's conjecture holds,
the Frobenius acts semi-simply at the eigenvalue $1$ on 
$H^{2n}_\et(\bar X,\Q_l(n))$, and
rational and numerical equivalence agree up to torsion on 
codimension $n$ cycles, for all $X$ and $n$.
\end{proposition}

\proof 
By the Weil conjectures and Gabber's theorem, the groups
$H^i_\et(X,\Q/\Z(n))$ and $\prod_ lH^i_\et(X,\Z_l(n))$
are finite unless $i=2n, 2n+1$. Thus the conjecture follows for 
$i\geq 2n+2$, because $H^{i-1}_\et(X,\Q/\Z(n))$ surjects onto the torsion group
$H^i_\et(X,\Z(n))$, which implies that they are of cofinite type, 
and finite for $i>2n+2$.

For $i<2n$, the finite group 
$H^{i-1}_\et(X,\Q/\Z(n))$ surjects onto $\Tor H^i_\et(X,\Z(n))$, and 
since $H^i_\et(X,\Z(n))\otimes\Q/\Z=0$, 
the short exact sequence 
$$ 0\to \Tor H^i_\et(X,\Z(n))\to  H^i_\et(X,\Z(n))\to  H^i_\et(X,\Q(n))\to 0$$
shows that $H^i_\et(X,\Z(n))$ modulo torsion is uniquely divisible. The vanishing
of this uniquely divisible subgroup is a restatement of Parshin's conjecture.

For $i=2n$ consider the sequence
$$ 0\to H^{2n}_\et(X,\Z(n))^{\wedge l}\to  H^{2n}_\et(X,\Z_l(n))
\to T_l H^{2n+1}_\et(X,\Z(n))\to 0.$$
The group of cofinite type  $H^{2n+1}_\et(X,\Z(n))$ is finite if 
and only if its Tate module vanishes if and only if the injection
$H^{2n}_\et(X,\Z(n))^{\wedge l}\to  H^{2n}_\et(X,\Z_l(n))$ of finitely
generated $\Z_l$-modules is an isomorphism. This is equivalent to
Tate's conjecture, because in the composition
$$CH^n(X)\otimes \Z_l\to  H^{2n}_\et(X,\Z(n))\otimes \Z_l \to 
H^{2n}_\et(X,\Z(n))^{\wedge l}$$
the first map is rationally an isomorphism, and the second map 
is surjective by Nakayama's Lemma. 

The final statement is \cite[Thm. 8.4]{ichweil}. Indeed, under all
those conjectures, Weil-\'etale cohomology and \'etale cohomology 
agree in degrees $\leq 2n$ and 
$H^{2n+1}_\et(X,\Z(n))\cong \Tor H^{2n+1}_W(X,\Z(n))$
by loc.cit. Thm. 7.1.
\endproof

\subsection*{Local fields}

\begin{theorem}
Let $k$ be a $p$-adic field of residue characteristic $p$. 
Then $H^i_\et(X,\Z(n))$ is
the direct sum of a finite group and a group which is uniquely $l$-divisible 
for all $l\not=p$, if one of the following two conditions holds
\begin{enumerate}
\item $X$ has good reduction and $i\not\in \{2n-1, 2n,2n+1, 2n+2\}$
\item $i\not \in \{n,\ldots , n+d+2\}$
\end{enumerate}
For $i=2n-1$ and $X$ with good reduction, or for $i=n$ in general, 
$\Tor H^i_\et(X,\Z(n)) $ is finite plus $p$-torsion. 
\end{theorem}

\proof
Assume first that $X$ has a smooth model $\mathcal X$ with closed fiber $Y$.
Then after inverting $p$, purity \cite{fujiwara} gives an exact sequence
$$\cdots  \to H^i_\et(\mathcal X,\Q/\Z(n))\to H^i_\et(X,\Q/\Z(n))\to 
H^{i-1}_\et(Y,\Q/\Z(n-1))\to\cdots . $$
The left hand term is isomorphic to $H^i_\et(Y,\Q/\Z(n))$
by the proper base change theorem. From the Weil conjectures
and Gabber's theorem \cite{gabber},
we know that $H^j_\et(Y,\Q/\Z(u))$ is finite for $j\not=2u,2u+1$,
and it follows that $H^i_\et(X,\Q/\Z(n))$ is finite 
for $i\not=2n-1,2n,2n+1$ after inverting $p$.
In particular, $H^i_\et(X,\Z(n))\otimes \Q/\Z$ is a $p$-group, 
and $\Tor H^{i+1}_\et(X,\Z(n))$ is the direct sum
of a $p$-group of cofinite type and a finite group, hence a direct summand \cite{kap}.
Finally, the sequence 
$$0\to \Tor H^{i}_\et(X,\Z(n))\to H^{i}_\et(X,\Z(n))\to 
H^{i}_\et(X,\Q(n))\to (p-\text{group})\to 0$$
show that $H^{i}_\et(X,\Z(n))$ modulo its torsion subgroup is uniquely $l$-divisible
for all $l\not=p$.

For general $X$, Kahn \cite[Thm.6 a)]{kahnfre} uses a strengthening of
the above argument to show that, for any $X$, 
$H^i_\et(X,\Q/\Z(n))$ is finite after localizing at $l\not=p$
for $i\not\in \{n,\ldots , n+d+1\}$. The rest of the proof continues
as above.
\endproof

\begin{example}[Weight $0,1$] 
The bounds in the good reduction case are sharp. If $X$ has a rational point, 
then for $r=[k:\Q_p]$, we obtain
\begin{align*}
i=2n-1 &\qquad  H^1_\et(X,\Z(1))\cong k^\times\cong \Z\oplus\Z_p^r\oplus (finite)\\
i=2n\phantom{ - 1 } &\qquad  H^0_\et(k,\Z(0))\cong \Z\\
&\qquad  
H^2_\et(X,\Z(1)) \cong \Pic(X)\cong (fin. gen.)\oplus \Z_p^{r\dim \Pic^0_X}\\
i=2n+1 &\qquad  H^3_\et(k,\Z(1))\cong \Br(k)\cong \Q/\Z\\
i=2n+2 &\qquad  H^2_\et(k,\Z(0))\cong 
\Hom(G_k,\Q/\Z)\cong \Q/\Z\oplus (\Q_p/\Z_p)^r\oplus (finite)
\end{align*}
The structure of the Picard group follows from a theorem of Mattuck.
The former groups have free summands, and the latter have
infinite torsion summands.

Moreover, a trace argument shows that $H^i_\et(k,\Z(n))$ is a direct summand 
of $H^{2d+i}_\et(X,\Z(d+n))$, i.e. 
the above phenomena occur for higher weights as well.
\end{example}

\begin{example}[Curves, weight $1$]
If $C$ is a curve of genus $g$ over a $p$-adic field, then
\begin{align*}
H^1_\et(C,\Z(1))&\cong k^\times\cong \Z\oplus\Z_p^r\oplus (finite)\\
H^2_\et(C,\Z(1))&\cong \Pic C\cong \Z\oplus \Z_p^{gr}\oplus (finite)\\
H^3_\et(C,\Z(1))&\cong \Br C\cong \Q/\Z \oplus (\Q_p/\Z_p)^{gr}\oplus (finite)
\end{align*}
The Brauer group is dual to the Picard group by Lichtenbaum
\cite{lichtenbaumcurve}.
\end{example}

\begin{example}[Weight $d+1$]\label{exw}
From the isomorphism
$$(\Tor H^2_\et(X,\Z))^* \cong H^1_\et(X,\Q/\Z)^* =\pi_1^{ab}(X),$$ 
and Poincar\'e duality $H^1_\et(X,\Z/m)^*\cong H^{2d+1}_\et(X,\Z/m(d))$, we
obtain a short exact sequence
$$ 0\to H^{2d+1}_\et(X,\Z(d+1))^\wedge \to \pi_1^{ab}(X) \to 
T H^{2d+2}_\et(X,\Z(d+1))\to 0.$$

If $X$ is a curve, then $H^3_\et(X,\Z(2))\cong SK_1(X)$, and 
we recover the exact sequence of Saito \cite[Thm. 2.6]{saitocurve}. 
In particular,  $\corank H^4_\et(X,\Z(2))=\rk H_1(\Gamma,\Z)= \rk H^1_\et(X,\Z)$,
which is zero if $X$ has good reduction. In arbitrary dimension, 
Yoshida shows \cite{yoshida} that 
$\corank H^{2d+2}_\et(X,\Z(d+1))$ is the dimension
of the maximal split torus of the Neron model of $Alb_X$ \end{example}

\begin{question}
%
1) If $X$ has good reduction and $i\leq 2n$, is $H^i_\et(X,\Z(n))$ 
the direct sum of a
finitely generated group and a finitely generated torsion free 
$\Z_p$-module (or at least a torsion free $\Z_{(p)}$-module)? 

2) If $X$ has good reduction and $i>2n$, is $H^i_\et(X,\Z(n))$ 
the direct sum of a cotorsion group
and a cofinitely generated torsion $\Z_{(p)}$-module?
\end{question}

Example \ref{exw} shows that the answer to the first question is negative if
$X$ does not have good reduction.

\begin{question}
Outside the range of the theorem,
the groups $H^i_\et(X,\Q_l/\Z_l(n))$ can have non-zero divisible
subgroups, and the groups 
$H^i_\et(X,\Q_p/\Z_p(n))$ have at least the same corank as
$H^{i-1}_\et(\bar X,\Q_p/\Z_p)$ for all $i$ and $n$ \cite[Cor. 11]{kahnfre}. 
Do these groups correspond to a non-torsion 
subgroup of $H^i_\et(X,\Z(n))$ or a torsion subgroup of $H^{i+1}_\et(X,\Z(n))$
in the short exact coefficient sequence
$$0\to H^i_\et(X,\Z(n))\otimes \Q_l/\Z_l \to H^i_\et(X,\Q_l/\Z_l(n))
\to H^{i+1}_\et(X,\Z(n))\{l\}\to 0 \quad ?$$
\end{question}

\begin{question} 
In the examples, the uniquely $l$-divisible groups appearing
have a structure of an $\mathcal O_k$-module. Is there 
such a structure in general on a subgroup (of countable index) of
$H^i_\et(X,\Z(n))$?
\end{question}

\subsection*{Arithmetic schemes}
Let $j:B\to C$ be an open subset of the spectrum of the ring of integers of a
number field or of a smooth and proper curve over a finite field, 
and let $S=C-B$ be the set of those places of $C$ not corresponding
to a point in $B$,
including the infinite places. Let $f:X\to B$ be smooth and proper, and
$d$ the dimension of $X$.
For a complex of sheaves $\f^\cdot$ on $X$ 
we define cohomology with compact support as the \'etale cohomology with compact
support on $B$ of $Rf_!\f^\cdot$, see \cite[\S 3]{kato} and \cite[II \S 2]{adt}. 
If $B$ is a curve over a finite field, or if we disregard $2$-torsion, 
then this is the cohomology of $j_!Rf_*\f^\cdot$ on $C$.
By construction, there is an exact sequence \cite[Prop. II 2.3]{adt}
\begin{equation}\label{cet}
\cdots \to H^i_c(X,\f^\cdot) \to H^i_\et(X,\f^\cdot) \to 
\bigoplus_{v\in S} H^i_\et(X\times_BK_v,\f^\cdot)
\to\cdots, 
\end{equation}
for $K_v$ the completion of the function field $K$ of $B$ at the place $v$;
in case of real fields, it is the Tate-modified cohomology of 
$R\Gamma_\et(X_{\mathbb C},\f^\cdot)$. For properties of Bloch's higher Chow groups
on smooth schemes over a Dedekind ring see \cite{ichdede}.
The following is the analog of Conjecture \ref{lconj}:

\begin{conjecture} (Lichtenbaum)
If $X$ is regular and proper over the ring of integers of a number field, 
then the groups  
$H^{i}_\et(X,\Z(n))$ 
are finitely generated for $i\leq 2n$,
finite for $i=2n+1$ and of cofinite type for $i\geq 2n+2$.
\end{conjecture}

Since the Tate-cohomology of $R\Gamma_\et(X_{\mathbb C},\Z(n))$ is easily seen
to be a finite $2$-group, this is equivalent to the same statement for 
$H^{i}_c(X,\Z(n))$. 
The following example shows that 
the analog statement does not hold for arbitrary $B$.

\begin{example}\label{dvrex}
Let $X=B$ be as above. Then we have the following isomorphisms and
exact sequences
\cite[II Prop. 2.1, 2.6, Cor. 2.11]{adt}:
\begin{enumerate}
\item[$n=0$:]
\begin{align*}
H^0_\et(B,\Z)&\cong \Z,\\
H^1_\et(B,\Z)&\cong 0,\\
H^2_\et(B,\Z)&\cong \Hom(\pi_1^{ab}(B),\Q/\Z).
\end{align*}
\item[$n=0$, compact support:]
$$ 0\to H^0_c(B,\Z)\to \Z \to \bigoplus_{v\in S} \Z \to H^1_c(B,\Z) \to 0,$$
$$ 0\to H^2_c(B,\Z)\to \pi_1^{ab}(B)^* \to \bigoplus_{v\in S} \pi_1^{ab}(K_v)^*
 \to H^3_c(B,\Z) \to H^3_\et(B,\Z) \to   0.$$
\item[$n=1$:]
\begin{align*}
H^1_\et(B,\Z(1))&\cong \Gamma(B,\G_m)^\times,\\
H^2_\et(B,\Z(1))&\cong \Pic(B),
\end{align*}
$$ 0\to H^3_\et(B,\Z(1))\to\bigoplus_{v\in S} \Br(K_v)
\to \Q/\Z \to H^4_\et(B,\Z(1))\to 0.$$
The last group vanishes if $B\not=C$; the higher groups are $2$-torsion.
\item[$n=1$, compact support:]
\begin{align*}
H^3_c(B,\Z(1))&\cong 0,\\
H^4_c(B,\Z(1))&\cong \Q/\Z,\\
H^i_c(B,\Z(1))&=0 ,\quad i>4.
\end{align*}
\end{enumerate}
\end{example}

The following is well-known, but we give a proof for the convenience 
of the reader:

\begin{proposition} If $n\geq 2$ and $B$ is as above, then 
$H^i_\et(B,\Z(n))$ is finitely generated. 

1) In the number field case,   
$H^1_\et(B,\Z(n))$ is finitely generated of rank
$r_1+r_2$ and $r_2$ if $n$ is odd and even, respectively,
and all other groups are finite.

2) In the function field case, the groups
$H^i_\et(B,\Z(n))$ are finite for all $i$.
\end{proposition}

\proof 
1) For $B=\Spec \O_F$, the case of (non-\'etale) motivic cohomology 
$H^i_\M(B,\Z(n))$ has been treated in \cite[Prop. 2.1]{KS}. By
the Rost-Voevodsky theorem and \cite{ichdede}, this gives the result
for $i\leq n+1$. 
Since the groups $H^i_\et(\F_q,\Z(n))$
are finite for $n\geq 1$ by Example \ref{ff}, the long exact sequence of 
\cite[Cor. 7.2]{ichduality} shows that the result holds for arbitrary $B$. 
If $i>n+1\geq 3$, the same long exact sequence shows that 
$$H^i_\et(B,\Z(n))\cong H^i_\et(K,\Z(n)) \cong H^{i-1}_\et(K,\Q/\Z(n))$$
which is (at most) a finite $2$-group.

2) Since $\Z/p(n)=0$ by \cite{ichmarc}, the groups are uniquely
$p$-divisible. Away from $p$, we can use 
the localization sequence and purity to conclude by the known
motivic cohomology of $C$ (Thm. \ref{ccc}) and finite fields 
(Example \ref{ff}).
\endproof

\begin{example} There are smooth and proper abelian schemes $X$ 
of dimension $3$ over $B$ such that $H^4_\et(X,\Z(2))$ contains a divisible
torsion group \cite{rossri}.
\end{example}

Lichtenbaum's conjecture suggests the following more general

\begin{conjecture}
The groups 
$H^{i}_\et(X,\Z(n))$ are a direct sum of a finitely generated group
and a group of cofinite type.
\end{conjecture}

The analog statement does not hold for $H^i_c(X,\Z(n))$ if $S$ contains a
finite prime $v$; for example the uncountable group $K_v^\times / K^\times$
maps onto $H^2_c(B,\Z(1))$.  
 
\bigskip
 
To end the paper, we consider an analog of the Tate-Shafarevich group:

\begin{question} Is the group
$$\Sh^{i,n}(X)= \ker  H^i_\et(X,\Z(n)) \to 
\bigoplus_{v\in S} H^i(X\times_BK_v,\Z(n))$$
finite if $S$ contains at least one finite prime? 
\end{question}

If $S$ does not contain any finite prime, then 
$\Sh^{4,1}(B)\cong \Q/\Z$, hence the condition is necessary.

\begin{example}
It is easy to see that 
$$\Sh^{0,0}(X)=\Sh^{1,0}(X)=\Sh^{0,1}(X)=\Sh^{1,1}(X)=0.$$
\end{example}


\begin{theorem}
If $f:X\to B$ is smooth and proper,
then $\Sh^{3,1}(X)$ is finite if and only if $\Sh(\Pic^{0,red}_{X_\eta/\eta})$ is finite.
\end{theorem}



\proof 
Recall that $\Z(1)\cong\G_m[-1]$.
Replacing $B$ by $f_*\mathcal O(X)$, we can assume that 
$f_*\G_m \cong \G_m$. Indeed, since $X$ is normal,
so is $f_*\mathcal O(X)$, hence is an open subset of the ring of integers
of a number field or smooth and proper curve over a  finite field again.
Moreover, $f': X\to \Spec f_*\mathcal O(X)$ is again smooth 
by SGA I, II Cor. 2.2, and proper. 
%
Consider the spectral sequence
$$ H^s_\et(B,R^tf_*\G_m)\Rightarrow H^{s+t}_\et(X,\G_m),$$
and compare it with the local situation:
$$ \begin{CD}
@>d_2>> E^{2,0}=H^2_\et(B,\G_m )=\Br(B) @>>> \Br(X)=H^2_\et(X,\G_m)@>\alpha>>\\
@.@VVV @V\beta VV \\
@>d_2>> \bigoplus_v H^2_\et(K_v,\G_m ) @>>> \bigoplus H^2_\et(X\times_BK_v,\G_m)
\end{CD}$$
The left vertical map is injective by Example \ref{dvrex}, and
the kernel of the lower horizontal is finite. Indeed, for each $v$,
$X\times_BK_v$ has a rational point over a finite extension of degree $d$,
and then the $l$-component of the kernel is cyclic of order dividing $d$.
We conclude that $\ker \alpha\cap \ker\beta$ is finite, and 
so it suffices to consider the quotient $\Br_0(X)= \Br(X)/\im \Br(B)$. 
From the spectral sequence we obtain a short exact sequence
$$ 0\to H^1_\et(B,R^1f_*\G_m) \to \Br_0(X)\to H^0_\et(B,R^2f_*\G_m) 
\stackrel{d_2}{\to }\cdots .$$
Let $g:\eta\to B$ be the inclusion of the generic point of $B$. 
Since $\Br(Y)\subseteq \Br(k(Y))$
for any regular integral scheme $Y$, we obtain that  
$R^2f_*\G_m\subseteq g_*R^2f_*\G_m|_\eta$. 
Hence 
$$H^0_\et(B,R^2f_*\G_m)\subseteq H^0_\et(B,g_*R^2f_*\G_m|_\eta)=
\Br(X\times_B\bar K)^{Gal(K)},$$
and the latter is a subgroup of $H^0_\et(K_v,R^2f_*\G_m)
\cong \Br(X\times_B\bar K_v)^{Gal(K_v)}$ for any $v$, because
$\Br(X_{\bar K})\subseteq  \Br(X_{\bar L})$ for an extension 
$\bar K\subseteq \bar L$ of 
algebraically closed fields by Proposition \ref{suslin}.

It remains to consider $H^1_\et(B,P)$ for $P=R^1f_*\G_m$. 
We claim that there is isomorphism $P\cong g_*g^*P$. Indeed, the map of 
stalks at a closed point of $B$
is $H^1_\et(X\times_B B_v^{sh},\G_m)\to H^1_\et(X\times_B K_v^{sh},\G_m)$,
for $B_v^{sh}$ the strict henselization of 
$B$ at $v$ and $K_v^{sh}$ its fields of fractions.
This map is an isomorphism by the localization sequence of higher
Chow groups because $X$ is smooth. 

Now consider the reduced connected component $A=\Pic^{0,red}_{X_\eta/\eta}$
of the Picard scheme $g^*P=\Pic_{X_\eta/\eta}$ at the generic point.
Then $g_*A$ is represented by the Neron model $\mathcal A $ of $A$.
Hence the exact sequence $0\to A\to g^*P \to \NS(X_\eta)\to 0$
induces an 
exact sequence $0\to \mathcal A\to g_*g^*P \to N\to 0$, where $N$ is a
subgroup of $g_* \NS(X_\eta)$. The resulting
long exact sequence of cohomology groups shows that 
$H^1_\et(B, \mathcal A)\to H^1_\et(B, g_*g^*P)\cong H^1_\et(B, P)$ differ
by a finite group. Similarly, $H^1_\et(K_v, A)\to H^1_\et(K_v, g^*P)$
differ by a finite group. Thus we can conclude by the isomorphism 
$$\Sh(A)\cong \ker H^1_\et(B, \A)\to \bigoplus_{v\in S} H^1_\et(K_v,A) 
$$
of \cite[II Lemma 5.5]{adt}.
\endproof

\small

{\sc  Department of Mathematics, Rikkyo University,\\
Nishi-ikebukuro, Toshimaku, Tokyo, Japan}

\textit{E-mail address:} {\tt geisser@rikkyo.ac.jp}


\begin{thebibliography}{SGA3}
\bibitem{blochesnault} {\sc S.\ Bloch, H.\ Esnault}, 
The coniveau filtration and non-divisibility for algebraic cycles. 
Math. Ann. 304 (1996), no. 2, 303--314. 

\bibitem{ctras} {\sc J.L.\ Colliot-Th\'el\`ene, W.\ Raskind}, 
$K_2$-cohomology and the second Chow group. 
Math. Ann. 270 (1985), no. 2, 165--199. 


\bibitem{sga4.5} {\sc P.\ Deligne}, Cohomologie \'etale. S\'eminaire de 
G\'eom\'etrie Alg\'ebrique du Bois-Marie SGA 412. 
Avec la collaboration de J. F. Boutot, A. Grothendieck, L. Illusie et 
J. L. Verdier. Lecture Notes in Mathematics, Vol. 569. Springer-Verlag, 
Berlin-New York, 1977. iv+312pp. 

\bibitem{fujiwara} {\sc K.\ Fujiwara}, A Proof of the Absolute Purity 
Conjecture (after Gabber), Advanced Studies in Pure Mathematics 36, 2002
Algebraic Geometry 2000, Azumino pp. 153--183.

\bibitem{gabber} {\sc O.\ Gabber}, Sur la torsion dans la cohomologie 
$l$-adique d'une vari\'et\`e. 
C. R. Acad. Sci. Paris S\'er. I Math. 297 (1983), no. 3, 179--182. 

\bibitem{ichtate} {\sc T.\ Geisser}, Tate's conjecture, algebraic cycles and 
rational $K$-theory 
in characteristic $p$. K-Theory 13 (1998), no. 2, 109--122. 

\bibitem{ichdede} {\sc T.\ Geisser}, Motivic cohomology over Dedekind rings. 
Math. Z. 248 (2004), no. 4, 773--794. 

\bibitem{ichweil} {\sc T.\ Geisser}, Weil-etale cohomology, 
Math. Ann. 330 (2004), 665--692.


\bibitem{ichparshin} {\sc T.\ Geisser}, 
Parshin's conjecture revisited. K-theory and noncommutative geometry, 
413--425, EMS Ser. Congr. Rep., Eur. Math. Soc., Z\"urich, (2008).

\bibitem{ichduality} {\sc T.\ Geisser}, Duality via cycle complexes,
Ann.\  of Math.\  (2) 172 (2010), no.\  2, 1095--1126.

\bibitem{ichmarc} {\sc T.\ Geisser, M.\ Levine}, 
The K-theory of fields in characteristic $p$. 
Invent. Math. 139 (2000), no. 3, 459--493.

\bibitem{grossuwa} {\sc M.\ Gros, N.\ Suwa} Application d'Abel-Jacobi 
$p$-adique et cycles alg\'ebriques. Duke Math. J. 57 (1988), no. 2, 579--613. 


\bibitem{kahnglr} {\sc B.\ Kahn}, The Geisser-Levine method revisited and 
algebraic cycles over a finite field. Math. Ann. 324 (2002), no. 3, 581--617. 

\bibitem{kahnfre} {\sc B.\ Kahn}, Some finiteness results for \'tale cohomology. 
J. Number Theory 99 (2003), no. 1, 57--73.

\bibitem{kap}{\sc I.\ Kaplansky}, Infinite abelian groups. 
Revised edition The University of Michigan Press, Ann Arbor, Mich. 1969 vii+95 pp.

\bibitem{kato}{\sc K.\ Kato}, A Hasse principle for two dimensional global
fields. Journal f\"ur die reine und angewandte Mathematik 366, 142--180.

\bibitem{korita} {\sc T.\ Kohrita}, Thesis, Nagoya University

\bibitem{KS} {\sc M.\ Kolster, J.W.\ Sands}, 
Annihilation of motivic cohomology groups in cyclic 2-extensions. 
Ann. Sci. Math. Qu\'ebec 32 (2008), no. 2, 175--187. 

\bibitem{lichtenbaumcurve} {\sc S.\ Lichtenbaum}, Duality theorems for curves 
over p-adic fields. Invent. Math. 7 (1969) 120--136. 

\bibitem{lichtenbaumMC} {\sc S.\ Lichtenbaum}, Values of zeta-functions at 
nonnegative integers. 
Number theory, Noordwijkerhout 1983 (Noordwijkerhout, 1983), 127--138, 
Lecture Notes in Math., 1068, Springer, Berlin, 1984.

\bibitem{lichtenbaumweil} {\sc S.\ Lichtenbuam}, The Weil-\'etale topology on schemes over finite fields, Compositio Math. 141 (2005) 689--702.

\bibitem{milneppart} {\sc J.\ Milne}, 
Values of zeta functions of varieties over finite fields. 
Amer. J. Math. 108 (1986), no. 2, 297--360.

\bibitem{adt} {\sc J.\ Milne}, Arithmetic duality theorems. Second edition. BookSurge, LLC, Charleston, SC, 2006. viii+339 pp. ISBN: 1-4196-4274-X

\bibitem{murrebar} {\sc J.\ Murre},  
Applications of algebraic K-theory to the theory of algebraic cycles. 
Algebraic geometry, Sitges (Barcelona), 1983, 216--261, 
Lecture Notes in Math., 1124, Springer, Berlin, 1985. 
 
\bibitem{murretor} {\sc J.\ Murre}, 
Algebraic cycles and algebraic aspects of cohomology and K-theory. 
Algebraic cycles and Hodge theory (Torino, 1993), 93--152, 
Lecture Notes in Math., 1594, Springer, Berlin, 1994.
 
\bibitem{rossri} {\sc A.\ Rosenschoen, V.\ Srinivas},
Torsion in the Lichtenbaum Chow group of arithmetic schemes.
Preprint.
 
\bibitem{saito} {\sc H.\ Saito}, 
Abelian varieties attached to cycles of intermediate dimension. 
Nagoya Math. J. 75 (1979), 95--119.

\bibitem{saitocurve} {\sc S.\ Saito}, Class field theory for curves over local fields. J. Number Theory 21 (1985), no. 1, 44--80. 

\bibitem{schoenma} {\sc C.\ Schoen}, Some examples of torsion in the 
Griffiths group. Math. Ann. 293 (1992), no. 4, 651--679. 

\bibitem{schoenjams} {\sc C.\ Schoen}, 
On the image of the l-adic Abel-Jacobi map for a variety over the 
algebraic closure of a finite field. 
J. Amer. Math. Soc. 12 (1999), no. 3, 795--838. 

\bibitem{schoenmrl} {\sc C.\ Schoen}, On certain exterior product maps of Chow groups. 
Math. Res. Lett. 7 (2000), no. 2-3, 177--194.

\bibitem{schoenjag} {\sc C.\ Schoen},
Complex varieties for which the Chow group mod $n$ is not finite. 
J. Algebraic Geom. 11 (2002), no. 1, 41--100. 


\bibitem{serre}{\sc J.P.\ Serre}, Morphisme universels et varieties d'albanese,
 Seminaire Chevalley,
1958-1959, Expose 10.

\bibitem{soulerec} {\sc C.\ Soul\'e}, The rank of \'etale cohomology of 
varieties over $p$-adic or number fields. 
Compositio Math. 53 (1984), no. 1, 113--131.

\bibitem{soulecurves} {\sc C.\ Soul\'e},  Groupes de Chow et K-th\'eorie de 
vari\'et\'s sur un corps fini. Math. Ann. 268 (1984), no. 3, 317--345.

\bibitem{yoshida} {\sc T.\ Yoshida}, Finiteness theorems in the class field 
theory of varieties over local fields. 
J. Number Theory 101 (2003), no. 1, 138--150.
\end{thebibliography}
\end{document}